\documentclass[a4paper,11pt]{article}
\usepackage[active]{srcltx}
\usepackage{array}
\usepackage{theorem}
\usepackage{amsmath,amscd,amssymb}
\usepackage{latexsym}
\usepackage{epsfig}
\usepackage{xypic}
\theorembodyfont{\sl}

\newtheorem{lemma}{Lemma}[section]

\newtheorem{proposition}[lemma]{Proposition}
\newtheorem{theorem}[lemma]{Theorem}
\newtheorem{corollary}[lemma]{Corollary}

\newtheorem{definition}[lemma]{Definition}

\newcommand{\CC}{\mathbb C}

\newcommand{\PP}{\mathbb P}
\newcommand{\QQ}{\mathbb Q}
\newcommand{\RR}{\mathbb R}

\newcommand{\ZZ}{\mathbb Z}

\renewcommand{\cD}{\mathcal D}

\newcommand{\cF}{\mathcal F}

\newcommand{\cI}{\mathcal I}

\renewcommand{\cL}{\mathcal L}
\newcommand{\cM}{\mathcal M}

\newcommand{\cO}{\mathcal O}

\newcommand{\isoto}{\stackrel{\sim}{\to}}

\renewcommand{\Tilde}{\widetilde}
\renewcommand{\Hat}{\widehat}

\newcommand{\Orth}{\mathop{\null\mathrm {O}}\nolimits}

\newcommand{\Orthtd}{\mathop{\null\mathrm{\widetilde{O}}}\nolimits}

\newcommand{\Opl}{{\Orth^+}}
\newcommand{\Og}{{\Orth_G}}
\newcommand{\TOpl}{{\Orthtd^+}}
\newcommand{\TOrth}{{\Orthtd}}

\newcommand{\lcm}{\mathop{\mathrm {lcm}}\nolimits}

\newcommand{\Aut}{\mathop{\mathrm {Aut}}\nolimits}

\newcommand{\id}{\mathop{\mathrm {id}}\nolimits}

\newcommand{\oriented}{\mathop{\null\mathrm {or}}\nolimits}

\newcommand{\Mon}{\mathop{\null\mathrm {Mon}}\nolimits}

\newcommand{\LA}{{L_A}}
\newcommand{\latt}[1]{{\langle{#1}\rangle}}
\renewcommand{\div}{\mathop{\mathrm {div}}\nolimits}
\newcommand{\Kthree}{\mathop{\mathrm {K3}}\nolimits}

\newcommand{\qedsymbol}{\mbox{$\Box$}}
\newcommand{\qed}{\unskip\nobreak\hfil\penalty50\hskip1em\hbox{}\nobreak
\hfill\qedsymbol\parfillskip=0pt\finalhyphendemerits=0}
\newenvironment{proof}{\begin{ProofwCaption}{Proof}}{\end{ProofwCaption}}
\newenvironment{ProofwCaption}[1]
 {\addvspace\theorempreskipamount \noindent{\it #1.}\rm}
 {\qed \par \addvspace\theorempostskipamount}

\begin{document}

\title{Moduli spaces of polarised symplectic O'Grady varieties and
 Borcherds products}
\author{V.~Gritsenko, K.~Hulek and G.K.~Sankaran}
\date{}
\maketitle
\begin{abstract}
We study moduli spaces of O'Grady's ten-dimensional irreducible
symplectic manifolds. These moduli spaces are covers of modular
varieties of dimension $21$, namely quotients of hermitian symmetric
domains by a suitable arithmetic group. The interesting and new
aspect of this case is that the group in question is strictly bigger
than the stable orthogonal group. This makes it different from both
the $\Kthree$ and the $\Kthree^{[n]}$ case, which are of dimension
$19$ and $20$ respectively.
\end{abstract}

\bigskip

\section{Introduction}
\label{sec0}

Irreducible symplectic manifolds are simply connected compact K\"ahler
mani\-folds which have a (up to scalar) unique $2$-form, which is
non-degenerate. In dimension two these are the $\Kthree$ surfaces. In
higher dimension there are, so far, four known classes of
examples. These are deformations of degree $n$ Hilbert schemes of
$\Kthree$ surfaces (the $\Kthree^{[n]}$ case), deformations of
generalised Kummer varieties, and two examples of dimensions $6$ and
$10$ due to O'Grady (\cite{OG2}, \cite{OG1}).

From the point of view of the Beauville lattice these examples fall
into two series. The first consists of $\Kthree$ surfaces, the
$\Kthree^{[n]}$ case and O'Grady's example of dimension $10$. The
Beauville lattices are the unimodular $\Kthree$-lattice
$L_{\Kthree}=3U \oplus 2E_8(-1)$, the lattice $L_{\Kthree} \oplus
\latt{-2(n-1)}$ and $L_{\Kthree} \oplus A_2(-1)$. The moduli spaces
of polarised irreducible symplectic manifolds of these classes are of
dimensions $19$, $20$ and $21$. The second series consists of
generalised Kummer varieties and O'Grady's $6$-dimensional variety
with Beauville lattices $3U \oplus \latt{-2 }$ and $3U \oplus
\latt{-2} \oplus \latt{-2} $ respectively. Here the dimensions of the
moduli spaces of polarised varieties are $4$ and $5$.

In order to describe moduli spaces of irreducible symplectic manifolds
one must first classify the possible types of the polarisation. We do
this in Section~\ref{sec3} for O'Grady's $10$-dimensional example. As
in the $\Kthree^{[n]}$ case we find that we have a split and a
non-split type. In this paper we shall mostly concentrate on the
split case, when the modular group is maximal possible, but we shall
also comment on the low degree non-split cases.

In the non-split case we expect Kodaira dimension $- \infty$ for the
three cases of lowest Beauville degree, namely $2d=12,\ 30,\ 48$. For
the next case of Beauville degree $2d= 66$ we prove general type: see
Corollary~\ref{prop-gt22}. The arguments used also suggest that
$2d=12$, $30$, $48$ might be the only degrees of non-split
polarisations giving unirational moduli spaces.

We should like to comment that there is a natural series consisting of
moduli of $\Kthree$ surfaces of degree $2$ (double planes branched
along a sextic curve), the non-split $\Kthree^{[2]}$ case of Beauville
degree $2d= 6$ (corresponding to cubic fourfolds and treated by Voisin
in \cite{Vo}) and O'Grady's example of dimension~$10$ with a non-split
polarisation of degree $12$. The lattices which are orthogonal to the
polarisation vector in this series are $2U \oplus 2E_8(-1) \oplus
A_n(-1)$ for $n=1,\ 2,\ 3$. It would be very interesting to find a
projective geometric realisation of O'Grady's $10$-dimensional
irreducible symplectic manifolds with non-split Beauville degree $12$.

In the split case we prove that the modular variety is of general type
for most degrees using the method of constructing low weight cusp
forms, as in the case of $\Kthree$ surfaces. The existence of such a
modular form proves that the modular variety is of general type,
provided the form vanishes along the branch divisors. We construct
these modular forms by using quasi-pullbacks of Borcherds' form
$\Phi_{12}$. There is, however, one important difference between the
split case for O'Grady varieties and the previous cases of $\Kthree$
surfaces \cite{GHS1} and the irreducible symplectic manifolds of
$\Kthree^{[n]}$-type~\cite{GHS2}. The modular group is now no longer
a subgroup of the stable orthogonal group: in fact it is a degree $2$
extension related to the root system $G_2$ (see
Theorem~\ref{thm-O(L,h)} and \eqref{eq-Og} below). This fact changes
considerably the geometry of the corresponding modular varieties. It
makes the case of the O'Grady varieties with a split polarisation very
interesting. We modify the original method of \cite{GHS1} and
\cite{GHS2} by considering involutions of the Dynkin diagrams and use
this to prove results for the split polarisation case (Sections
\ref{sec4}--\ref{sec5}). Here we make strong use of the
classification of lattices of small rank and determinant (see
Conway-Sloane \cite{CS}).

The case of Beauville degree $2d=2^n$ is exceptional because of very
special relations between the root systems $E_6$ and $F_4$. We cannot
obtain any results about the birational type of these modular
varieties. However, if we take the double cover given by the stable
orthogonal group, we can prove general type with the only exceptions
the split polarisations $2d=2$, $4$, $8$.

The geometry of roots is very special in this case and quite different
from the $\Kthree$ and the $\Kthree^{[n]}$ case. Because of some very
special coincidences we require no explicit Siegel type formulae for
the representation of an integer by a lattice, nor do we have to
enlist the help of a computer.
\smallskip

{\noindent {\bfseries Acknowledgements}:
We should like to thank Eyal Markman for informative
conversations on monodromy groups. We are grateful for financial support
under grants DFG Hu/337-6 and ANR-09-BLAN-0104-01.
The authors would like to thank Max-Planck-Institut
f\"ur Mathematik in Bonn for support and for providing excellent
 working conditions. \par}

\section{Irreducible symplectic manifolds and moduli}
\label{sec1}

We first recall the following.
\begin{definition}\label{def:irrsymplectic}
 A complex manifold $X$ is called an {\em irreducible symplectic
 manifold} or {\em hyperk\"ahler manifold} if the following
 conditions are fulfilled:
\begin{itemize}
\item[{\rm (i)}] $X$ is a compact K\"ahler manifold;
\item[{\rm (ii)}] $X$ is simply-connected;
\item[{\rm (iii)}] $H^0(X,\Omega^2_X) \cong \CC \sigma$ where $\sigma$
 is an everywhere nondegenerate holomorphic $2$-form.
\end{itemize}
\end{definition}
It follows from the definition that $X$ has even complex dimension,
$\dim_{\CC}(X)=2n$, and that the canonical bundle $\omega_X$ is
trivial (a trivializing section is given by $\sigma^n$). Moreover, the
irregularity $q(X)=h^1(X,\cO_X)=0$. Irreducible symplectic manifolds
are, together with Calabi-Yau manifolds and abelian varieties, one of
the building blocks of compact K\"ahler manifolds with trivial
canonical bundle (complex Ricci flat manifolds). In dimension~$2$ the
irreducible symplectic manifolds are the $\Kthree$ surfaces. So far
only four deformation types of such manifolds have been found. These
are (deformations of) Hilbert schemes of points on $\Kthree$ surfaces
(also called irreducible symplectic manifolds of
$\Kthree^{[n]}$-type), (deformations of) generalised Kummer varieties
and two types of examples constructed by O'Grady (see \cite{OG1},
\cite{OG2}).

For a $\Kthree$ surface $S$ the intersection form defines a non-degenerate,
symmetric bilinear form on the second cohomology $H^2(S,\ZZ)$, giving this
cohomology group the structure of a lattice. More precisely
$$
H^2(S,\ZZ) \cong 3 U \oplus 2 E_8(-1) = L_{\Kthree}
$$
where $U$ is the hyperbolic plane and $E_8(-1)$ is the unique even,
negative definite unimodular lattice of rank $8$. Similarly, one can
also define a lattice structure on $H^2(X,\ZZ)$ for all irreducible
symplectic manifolds $X$, called the \emph{Beauville lattice}. The
easiest way to define this is the following. There exists a positive
constant $c$, the \emph{Fujiki constant}, such that the quadratic form
$q$ on $H^2(X,\ZZ)$ defined by $(\alpha)^{2n}=cq(\alpha)^n$ is the
quadratic form of a primitive non-degenerate symmetric bilinear
form. This form has signature $(3,b_2(X)-3)$.

Let $L$ be an abstract lattice isomorphic to the Beauville lattice of an
irreducible symplectic manifold. This defines a period domain
$$
\Omega=\{ [x] \in \PP(L \otimes \CC ) \mid (x,x)=0, \, (x,\bar{x}) >0 \}.
$$
Given a marking on an irreducible symplectic manifold, i.e.\ an
isometry $\phi\colon H^2(X,\ZZ) \isoto L$, one can define the period
point of $X$ as the point in $\Omega$ defined by the line
$\phi_{\CC}(H^{2,0}(X))$. As in the $\Kthree$ case, irreducible
symplectic manifolds are unobstructed and local Torelli holds: that
is, the period map of the Kuranishi family is a local isomorphism (see
\cite{Be}). Moreover Huybrechts \cite{Huy} proved surjectivity of the
period map.

We are interested in moduli of polarised irreducible symplectic
manifolds. By a polarisation we mean a primitive ample line bundle
$\cL$ on $X$ and we call $h=c_1(\cL) \in H^2(X,\ZZ)$ the polarisation
vector. Since $\cL$ is ample, the Beauville degree $q(h)$ is strictly
positive. Note that the geometric degree of the polarisation is
$cq(h)^n$.

In order to discuss moduli spaces of polarised irreducible symplectic
varieties, one has to fix discrete data. These are firstly the
Beauville lattice and the Fujiki invariant (which together determine
the so-called numerical type of an irreducible symplectic manifold)
and secondly the type of the polarisation. Since the Beauville lattice
$L$ of an irreducible symplectic manifold is, in general, not
unimodular, we cannot expect that any two polarisation vectors of the
same degree are equivalent under the orthogonal group $\Orth(L)$. (The
case of $\Kthree$ surfaces is an exception, since the
$\Kthree$-lattice is unimodular.) In general there will be several,
but finitely many, $\Orth(L)$-orbits of such vectors. We call the
choice of such an orbit the choice of a \emph{polarisation
 type}. Given a polarisation type we fix a representative $h \in L$
of it and consider the lattice $L_h = h^{\perp}_L$, which has
signature $(2,b_2(X)-3)$, and defines a homogeneous domain
$$
\Omega_h=\Omega(L_h)=\{ [x] \in \PP(L_h \otimes \CC ) \mid (x,x)=0,
\, (x,\bar{x}) >0 \}.
$$
This is a type IV bounded symmetric hermitian domain. It is of
dimension $b_2(X)-3$ and has two connected components
$$
\Omega(L_h)=\cD(L_h) \coprod \cD(L_h)'.
$$
The orthogonal group $\Orth(L_h)$ of the lattice $L_h$ has an
index~$2$ subgroup $\Orth^+(L_h)$ that fixes the components
$\cD(L_h)$ and $\cD(L_h)'$. We also need the group
\begin{equation}\label{eq-Oh}
\Orth(L,h) = \{ g\in \Orth(L) \mid g(h)=h \}.
\end{equation}
Since this group maps the orthogonal complement $L_h$ to itself, we
can consider it as a subgroup of $\Orth(L_h)$. Let $\Orth^+(L,h)=
\Orth(L,h) \cap \Orth^+(L_h)$.

Let $\cM_h$ be the moduli space of polarised irreducible symplectic
manifolds $(X,\cL)$ where $X$ has numerical data as chosen above and
where $\cL$ is a primitive ample line bundle such that $c_1(\cL)$ is
of the given polarisation type. This moduli space exists by Viehweg's
general theory as a quasi-projective variety. We do not know how many
components $\cM_h$ has, but Propostion~\ref{prop:modulicovering} below
allows us to work with each component separately.
\begin{proposition}\label{prop:modulicovering}
Every component $\cM_h^0$ of the moduli space $\cM_h$ admits a
dominant finite-to-one morphism
$$
\varphi\colon \cM^0_h \to \Orth^+(L,h) \backslash \cD(L_h).
$$
\end{proposition}
\begin{proof}
See \cite[Theorem 1.5]{GHS2}.
\end{proof}
This is the starting point of our investigations. The importance of
this result is that if the quotient $\Orth^+(L,h) \backslash \cD(L_h)$
is of general type, then so is $\cM^0_h$. We shall use this in
Sections~\ref{sec4} and \ref{sec5} to prove the main result of this
paper.

For some irreducible symplectic manifolds, such as irreducible
symplectic manifolds of $\Kthree^{[n]}$-type, the situation can be
improved by introducing the group $\Mon^2(X) \subset
\Orth(H^2(X,\ZZ))$, which is the group generated by the monodromy
group operators acting on the second cohomology. This group was
studied intensively by Markman (\cite{Mar1}, \cite{Mar2},
\cite{Mar3}). If it is a normal subgroup, then it defines a subgroup
$\Mon^2(L) \subset \Orth(L)$. One can then show (the proof of
\cite[Theorem 2.3]{GHS2} for the $\Kthree^{[n]}$-type goes through
unchanged) that one can factor the map $\varphi$ from
Proposition~\ref{prop:modulicovering} as follows:
\begin{equation}\label{equ:monodromyfactorisation}
\xymatrix{ {\cM^{0}_{h}} \ar[r]^(.2){\tilde\varphi} \ar[dr]^{\varphi} &
{(\Mon^2(L) \cap \Orth^+(L,h)) \backslash \cD(L_h)}
\ar[d]\\
& {\Orth^+(L,h) \backslash \cD(L_h).} }
\end{equation}

\section{O'Grady's $10$-dimensional example}
\label{sec2}

O'Grady constructed his $10$-dimensional irreducible symplectic
manifolds using moduli spaces of sheaves on $\Kthree$ surfaces. More
precisely, let $S$ be an algebraic $\Kthree$ surface and consider the
rank $2$ sheaves $\cF$ on $S$ with trivial first Chern class
$c_1(\cF)=0$ and second Chern class $c_2(\cF)=4$. Let $H$ be a
sufficiently general polarisation, i.e.\ a polarisation such that
there is no non-trivial divisor class $C$ with $C.H=0$ and $C^2 \geq
-4$. It is easy to find examples: in particular every projective
$\Kthree$ surface with Picard number~$1$ has such a polarisation. Let
$\cM_4$ be the moduli space of $H$-semistable sheaves. This is a
singular variety whose smooth part carries a symplectic structure. The
singularities occur at the semi-stable sheaves and these are sums of
ideal sheaves $\cI_Z \oplus \cI_W$ where $Z$ and $W$ are
$0$-dimensional subschemes of $S$ of length~$2$. O'Grady then
considers Kirwan's desingularisation $\Hat{\cM}_4$ which has a
canonical form vanishing on an irreducible divisor. He shows that this
divisor is a $\PP^2$-bundle whose normal bundle has degree $-1$ on
each $\PP^2$. Hence it can be contracted and the resulting $4$-fold
$\Tilde{\cM}_4$ is O'Grady's irreducible symplectic manifold of
dimension $10$. It has second Betti number $b_2=24$. This also shows
that these varieties have $22$ deformation parameters and hence there
are deformations of $\Tilde{\cM}_4$ which do not arise from
deformations of the underlying $\Kthree$ surface.

In the case of O'Grady's $10$-dimensional examples the Beauville
lattice is (as an abstract lattice) of the form:
$$
L=3U \oplus 2E_8(-1) \oplus A_2(-1)
$$ where $A_2(-1)$ is the negative definite root lattice associated to
$A_2$. The Fujiki invariant of O'Grady's $10$-dimensional example is
$c=945$. This was shown by Rapagnetta \cite{Ra}. Since the second
cohomology of $\Kthree$ surfaces is of the form $L=3U \oplus 2E_8(-1)$
and the Beauville lattices of irreducible symplectic manifolds of
$\Kthree^{[n]}$-type are of the form $L=3U \oplus 2E_8(-1) \oplus
\latt{-2(n-1)} $, one can see O'Grady's $10$-dimensional example as
the third type in a series. We previously treated the case of
$\Kthree$ surfaces in \cite{GHS1} and the case of polarised varieties
of $\Kthree^{[n]}$-type in \cite{GHS2}, where we restricted ourselves
to the case of split polarisations (see \cite[Example 3.8]{GHS1} for a
definition and details).

In the $10$-dimensional case the situation with respect to the
monodromy group is as follows. Let $\Orth^{\oriented}(L)$ be the group
of oriented orthogonal transformations of $L$ (see \cite[Section
  4.1]{Mar1} and in particular Remark~4.3 for a definition of oriented
orthogonal transformations). By a result of Markman (unpublished) it
is known that $\Mon^2(L)= \Orth^{\oriented}(L)$.

Since $\Orth(L,h) \cap\Orth^{\oriented}(L) = \Orth^+(L,h)$ the
factorisation~\eqref{equ:monodromyfactorisation} does not, unlike in
some cases of $\Kthree^{[n]}$-type, improve the situation. In view of
Verbitsky's results \cite{Ve} we conjecture that the map $\varphi \colon
\cM^0_h \to \Orth^+(L,h) \backslash \cD(L_h)$ from
Proposition~\ref{prop:modulicovering} is indeed an open embedding.

There are two differences between the cases treated previously and
this case. Firstly, the arithmetic group in question is no longer
necessarily a subgroup of the stable orthogonal group (see
Section~\ref{sec3}). Secondly, the discriminant group of the lattices
orthogonal to a polarisation vector is no longer cyclic. This
requires new considerations concerning the quasi-pullbacks of the
Borcherds form. We would also like to point out that the lattice
theoretic part of this case is very different from the previous
papers. The geometry of roots is very special here, and as a result we
need neither arguments from analytic number theory nor any kind of
Siegel formulae. The root geometry arguments in this paper are all
elementary, but they are far from trivial.

\section{The modular orthogonal group and the root system $G_2$}
\label{sec3}

In this section we determine the modular group associated to the
moduli spaces of polarised O'Grady varieties (see
Theorem~\ref{thm-O(L,h)} below). A polarisation corresponds to a
primitive vector $h$ with $h^2=2d>0$ in
\begin{equation}\label{eq-L}
\LA=3U\oplus 2E_8(-1)\oplus A_2(-1).
\end{equation}
For any even lattice $L$ we denote the discriminant group of $L$ by
$D(L)=L^{\vee}/L$ where $L^\vee$ is the dual lattice of $L$. The
discriminant group carries a discriminant quadratic form $q_L$ (if $L$
is even) with values in $\QQ/2\ZZ$. The orthogonal group of the finite
discriminant form is denoted by $\Orth(D(L))$. If $g\in \Orth(L)$ we
denote by $\bar g$ its image in $\Orth(D(L))$. The \emph{stable
 orthogonal group} $\Tilde\Orth(L)$ is defined by
$$
\Tilde{\Orth}(L) = \ker (\Orth (L) \to \Orth(D(L))).
$$
If $h\in L$ its \emph{divisor} $\div(h)$ is the positive generator
of the ideal $(h,L)\subset \ZZ$. Therefore $h^*=h/\div(h)$ is a
primitive element of the dual lattice $L^\vee$ and $\div(h)$ is a
divisor of $\det (L)$.

For the lattice $\LA$ of \eqref{eq-L}, $D(\LA)\cong
D(A_2(-1))=\latt{\bar c}$ is the cyclic group of order $3$ and
$q_\LA(\bar c)=\frac{2}{3}\mod 2\ZZ$.

For any $h\in \LA$ with $h^2>0$ and $L_h=h^\perp_\LA$ we determine the
structure of the modular group $\Opl(\LA,h)=\Orth(\LA,h)\cap
\Orth^+(L_h)$ (see \eqref{eq-Oh} and
\eqref{equ:monodromyfactorisation}).
We have $\det(\LA)=3$, so $\div(h)$ divides $(2d,3)$.

\begin{theorem}\label{thm-O(L,h)}
Let $h\in \LA$ be a primitive vector of length $h^2=2d>0$. The
orthogonal complement $L_h=h^\perp_\LA$ is of signature $(2,21)$. If
$\div(h)=3$ then
$$
L_h\cong L_{Q}=2U\oplus 2E_8(-1)\oplus Q(-1),
$$
where $Q(-1)$ is a negative definite even integral ternary
quadratic form of determinant $-2d/3$. Its discriminant group
$D(Q(-1))\cong D(L_h)$ is cyclic of order $2d/3$ and
$$
\Opl(\LA,h)\cong \TOpl(L_h).
$$
If $\div(h)=1$, then $L_h\cong L_{A, 2d}$ where
$$
L_{A, 2d}=2U\oplus 2E_8(-1)\oplus A_2(-1)\oplus \latt{-2d},
$$
$$
D(L_h)\cong D(A_2(-1))\oplus D(\latt{-2d}),
$$
and
$$
\Opl(\LA,h)\cong \Og(L_{A ,2d})=\{ g\in \Opl(L_{A,2d})\,|\ \bar
g|_{D(\latt{-2d})}=\id \}.
$$
Any totally isotropic subgroup of $D(A_2(-1))\oplus D(\latt{-2d})$ is cyclic.
\end{theorem}
\smallskip

A polarisation determined by a primitive vector $h_d$ with
$\div(h_d)=1$ is called \emph{split}. We note that if $(3,d)=1$ then
the polarisation is always split. If $3\vert d$ then the polarisation
$h=2d$ is split if and only if the discriminant group of $L_h$ is not
cyclic. In the split case the modular group $\Og(L_{A ,2d})$ is larger
than the stable orthogonal group $\TOpl(L_h)$ because the elements of
$\Og(L_{A ,2d})$ induce trivial action only on the second component of
the discriminant group $D(L_h)\cong D(A_2(-1))\oplus D(\latt{-2d})$.

We recall that
$$
[\Orth (A_2): W(A_2)]=2
$$
where $\Orth (A_2)$ is the orthogonal group of the lattice $A_2$ and
$W(A_2)$ is the Weyl group generated by reflections with respect to
the roots of $A_2$.  The group $O(A_2)$ contains also reflections with
respect to the vectors of square $6$. The $2$- and $6$-roots of the
lattice $A_2$ form together the root system $G_2$ and $\Orth
(A_2)=W(G_2)$ (see \cite{Bou}).

For any vector $l\in L_h$ with $l^2<0$ the reflection $\sigma_l$ with
respect to $l$ belongs to $\Orth^+(L_h\otimes \RR)$. In particular,
$O(A_2(-1))=W(G_2(-1))$ is a subgroup of $\Orth^+(L,h)$. Therefore
\begin{equation}\label{eq-Og}
\Og(L_{A ,2d})/\TOpl(L_h)\cong W(G_2(-1))/W(A_2(-1))\cong \ZZ/2\ZZ.
\end{equation}
We note that in the case of polarised $\Kthree$ surfaces or of
polarised symplectic manifolds of $\Kthree^{[n]}$-type the modular
group of the corresponding modular varieties is identical to a stable
orthogonal group (see \cite{GHS2}). The degree~$2$ extension of the
stable orthogonal group changes the geometry of the modular varieties
considerably.  This can be compared to the case of the moduli spaces
of $(1,p)$-polarised abelian and Kummer surfaces (see \cite{GH}).
\smallskip

Theorem \ref{thm-O(L,h)} shows the difference between split and
non-split polarisations. To prove it we study the orbits of vectors in
$L$. Using the standard discriminant group arguments (see \cite{Nik}
and the proof of Proposition 3.6 in \cite{GHS1}) we get
\begin{lemma}\label{lem-Lh}
Let $L$ be any non-degenerate even integral lattice and let $h\in \LA$
be a primitive vector with $h^2=2d>0$. If $L_h$ is the orthogonal
complement of $h$ in $\LA$ then
$$
\det L_h=\frac {(2d)\cdot \det \LA}{\div (h)^2}.
$$
\end{lemma}

A proof of the following classical result, known as the \emph{Eichler
  criterion}, is given in \cite[Proposition 3.3]{GHS4}.
\begin{lemma}\label{Eichler}
Let $L$ be a lattice containing two orthogonal isotropic planes. Then
the $\Tilde\Orth(L)$-orbit of a primitive vector $l\in L$ is
determined by two invariants: its length $l^2=(l,l)$, and its image
$l^*+L$ in the discriminant group $D(L)$.
\end{lemma}
\smallskip

According to this Lemma~\ref{Eichler}, all primitive $2d$-vectors
$l\in\LA$ with $\div(l)=1$ belong to the same $\TOrth(\LA)$-orbit. If
$\div(l)=3$ then $l^*+\LA$ is a generator of $D(\LA)=D(A_2(-1))$.
Therefore there are two $\TOrth(\LA)$-orbits of such vectors. An
element of $W(G_2(-1))$ makes these two $\TOrth(\LA)$-orbits into one
$\Orth(\LA)$-orbit.

\begin{lemma}\label{lem-nspl}
If $h_{2d}$ is a vector of a non-split polarisation then $2d\equiv 12
\mod 18$. For any positive even integer $2d$ satisfying this
congruence there exists a primitive $h_{2d}\in \LA$ with
$\div(h_{2d})=3$.
\end{lemma}
\begin{proof}
We put $h_{2d}=u+xa+yb\in \LA$, where $u\in 3U\oplus 2E_8(-1)$ and
$xa+yb\in A_2(-1)=\latt{a,b}$, where $a$, $b$ are simple roots of
$A_2(-1)$.  Any primitive vector of a unimodular lattice has divisor
$1$. Therefore $u=3v$ with $v\in 3U\oplus 2E_8(-1)$. A straightforward
calculation shows that $\div(xa+yb)$ is divisible by $3$ if and only
if $x+y\equiv 0\mod 3$.  We have $x\equiv \pm 1\mod 3$ and $y\equiv
\mp 1\mod 3$ since $h_{2d}$ is primitive. Therefore
$$
h_{2d}^2=9v^2-2(x+y)^2+6xy\equiv 12 \mod 18.
$$
To construct a polarisation vector of degree $18n-6$ we take a
vector $h=3nu_1+3u_2+(2a+b)$ where $U=\latt{u_1,u_2}$ is the first
hyperbolic plane in $\LA$.
\end{proof}

Now we can calculate $L_h$. If the polarisation is non-split we take
the vector $h_{2d}\in U\oplus A_2(-1)$ indicated above. We denote by
$Q(-1)$ the orthogonal complement of $h_{2d}$ in $U\oplus
A_2(-1)$. According to Lemma \ref{lem-Lh} it is an even integral
negative definite lattice of rank $3$ and of determinant $-2d/3$,
i.e.\
$$
L_h\cong 2U\oplus 2E_8(-1)\oplus Q(-1),\qquad \det Q(-1)=-\frac{2d}3.
$$
To prove that $D(Q)$ is cyclic we consider
$$
\latt{h} \oplus L_h\subset \LA \subset L_A^\vee \subset
\latt{\frac{1}{2d}h} \oplus L_h^\vee.
$$
The lattice $L$ defines the finite subgroup
$$
H=\LA/(\latt{h} \oplus L_h) < D(\latt{h})\oplus D(L_h).
$$
We have $|H|=\det L_h=2d/3$ because $H\cong (\latt{\frac{1}{2d}h}
\oplus L_h^\vee)/L_A^\vee$. The projections
\begin{alignat}{2}\label{eq-proj-pH}
p_h\colon H \to D(\latt{h}),\qquad & p_{L_h}\colon H\to D(L_h)
\end{alignat}
are injective because $\latt{h}$ and $L_h$ are primitive in $\LA$ (see
\cite[Prop. 1.5.1]{Nik}). Therefore $H\cong D(L_h)$ and $H$ is
isomorphic to a subgroup of the cyclic group $D(\latt{h})$.
\smallskip

To determine $\Orth(\LA,h)$ we consider the action of elements of this
group on the discriminant group.  Any $g\in \Orth(\LA,h)$ acts on
$\latt{h}^\vee\oplus L_h^\vee$ and induces an element $\bar g\in
\Orth(D(\LA))$. Moreover $\bar g$ acts on the subgroup $H$. For any
$\bar a \in p_h(H)$ there exists a unique $\bar b\in p_{L_h}(H)$ such
that $\bar a+\bar b\in H$. The action of $\bar g$ on $D(\latt{h})$ is
trivial.  Therefore it is also trivial on the second component $\bar
b\in p_{L_h}(H)$.  But $p_{L_h}(H)$ is isomorphic to the whole group
$D(L_h)$ if $\div(h)=3$.  Therefore $\Orth(\LA,h)\cong
\Tilde\Orth(\LA_h)$. This proves the statement of Theorem
\ref{thm-O(L,h)} in the non-split case.

For a split polarisation we can take $h_{2d}=du_1+u_2\in U$. Then
$(h_{2d})^\perp_{U}\cong \latt{-2d}$ and
$$
L_h\cong 2U\oplus 2E_8(-1)\oplus A_2(-1)\oplus \latt{-2d}.
$$
Then $|H|=2d$, $p_{L_h}(H)\cong D(\latt{-2d})$ and $\bar g$ acts
trivially on $D(\latt{-2d})$.

\smallskip

To finish the proof of Theorem \ref{thm-O(L,h)} we analyse the
isotropic elements of the discriminant group $D(A_2(-1))\oplus
D(\latt{-2d})$ of the lattice $L_h$ in the split case. If $(3,d)=1$,
then the latter group is cyclic. So we assume that $3|d$.

Let $\bar {l}=(\pm \bar c, \frac{x}{2d}\bar h)$ where $\bar c$ is a
generator of $D(A_2(-1))$ and $x$ is taken modulo $2d$.  We put
$d=3d_0=3ef^2$ where $e$ is square free.  It is easy to see that $\bar
l$ is isotropic if and only if $x=2yef$, where $y$ is taken
modulo~$3f$, and $1+ey^2\equiv 0$ mod~$3$.  The element $\bar l$ is
isotropic if and only if
$$
\frac{2}3+\frac{x^2}{2d}\equiv 0\mod 2d.
$$
Then
$$
4d+3x^2\equiv 0 \mod 12d \quad {\rm \ or\ }\quad 12ef^2+3x^2\equiv 0
\mod 36ef^2.
$$
We see that $x=2x_0$ and $ef^2+x_0^2\equiv 0$ mod~$3ef^2$.  Therefore
$x_0\equiv 0$ mod~$ef$ and $x=2x_0=2efy$ where $y$ is taken
modulo~$3f$ and
$$
1+ey^2\equiv 0\mod 3.
$$
The last congruence is true if and only if
$$
e\equiv 2\mod 3\quad{\rm and }\quad y\not\equiv 0\mod 3.
$$
We proved that for $d=3ef^2$ the isotropic elements with non trivial
first component are $(\pm \bar c,\ \frac{y}{3f}\bar h)$. All these
elements belong to the union of two totally isotropic cyclic groups
generated by $(\bar c,\ (\bar h/3f))$ and by $(\bar c,\ -(\bar
h/3f))$. If a subgroup of the discriminant group contains two
isotropic elements $(\bar c, y_i(\bar h/3f))$, where $y_1\not\equiv
y_2$ mod~$3$, then $(\bar 0, (y_1-y_2)(\bar h/3f))$ is not isotropic
because
$$
\frac{6ef^2(y_1-y_2)^2}{9f^2}=\frac {2e(y_1-y_2)^2}{3}\not\equiv 0\mod
2\ZZ.
$$
Thus Theorem \ref{thm-O(L,h)} is proved.
\smallskip

\noindent{\bf Example 1.} \emph{The smallest non-split polarisations
  $12$, $30$, $48$, $66$.} In the non-split case the isomorphism class
of the lattice $L_h$ with $h^2=2d$ is uniquely defined by the genus of
the ternary form $Q$ of determinant $2d/3$. For the small
polarisations of this example the genus of $Q$ contains only one
class. The corresponding classes can be found in \cite[Table
  I]{CS}. We give a modified description of them using the language of
root lattices, indicating the maximal root subsystem in the lattices
$Q$ and $Q^\perp_{E_8}$:
\begin{alignat*}{3}
\det Q=\ 4, &\quad Q=A_3,&\quad Q^\perp_{E_8}\cong
D_5,\qquad\qquad\qquad\quad\ \\
\det Q=10, &\quad Q=(A_1)^\perp_{A_4},
&\quad Q^\perp_{E_8}\cong A_1\oplus A_4,\qquad\qquad\quad\\
\det Q=16, &\quad Q \supset A_2\oplus \latt{48},
&\quad Q^\perp_{E_8}\supset A_4\oplus \latt{48},\qquad\qquad\ \\
\det Q=22, &\quad Q\supset A_2\oplus \latt{66}, & \quad Q^\perp_{E_8}\supset
A_3\oplus A_1\oplus \latt{44}.\qquad
\end{alignat*}

\section{Cusp forms of small weight and the Borcherds form $\Phi_{12}$}
\label{sec4}

Now we can formulate the main theorem of the paper.

\begin{theorem}\label{thm-main}
Let $d$ be a positive integer not equal to $2^n$ with $n\ge 0$.
Then the modular variety
$$
M_{A,2d}=\Og(L_{A ,2d})\setminus \cD(L_{A, 2d})
$$
is of general type.  Every component $\cM_h^0$ of the moduli space
$\cM_h$ of ten-dimensional polarised O'Grady varieties with split
polarisation $h$ of Beauville degree $h^2=2d\ne 2^{n+1}$ is of general
type.
\end{theorem}

\noindent
{\bf Remark.} In Corollary \ref{prop-gt22} below we prove general type
of the moduli spaces $\cM_h^0$ for the fourth non-split polarisation,
of Beauville degree $66$ (see Example 1 of \S\ref{sec3}).
\smallskip

According to Proposition \ref{prop:modulicovering} it is enough to
prove the main Theorem \ref{thm-main} for the modular varieties
$$
M_{A, 2d}=\Og(L_{A ,2d})\setminus \cD(L_{A, 2d})\quad {\rm or }\quad
M_{Q}^{(2d)}=\TOpl(L_{Q})\setminus \cD(L_{Q})
$$
(see notations of Theorem \ref{thm-O(L,h)}). The dimension of the
modular variety $M_{A, 2d}$ is $21$, which is larger than
$8$. Therefore we can use the \emph{low weight cusp form trick} from
\cite{GHS1}.

Let $L$ be an even integral lattice of signature $(2,n)$ with $n\ge
3$. A modular form of weight $k$ and character $\det$ with respect to
a subgroup $\Gamma<\Orth^+(L)$ of finite index is a holomorphic
function $F\colon\cD(L)^\bullet\to \CC$ on the affine cone
$\cD(L)^\bullet$ over $\cD(L)$ such that
\begin{equation*}\label{mod-form}
 F(tZ)=t^{-k}F(Z)\quad \forall\,t\in \CC^*\text{\qquad and\qquad}
 F(gZ)=\det(g)F(Z)\quad \forall\,g\in \Gamma.
\end{equation*}
A modular form is a cusp form if it vanishes at every cusp. Cusp forms
of character $\det$ vanish to integral order at any cusp (see
\cite{GHS4}). We denote the linear spaces of modular and cusp forms of
weight~$k$ and character $\det$ for $\Gamma$ by $M_k(\Gamma,\det)$ and
$S_k(\Gamma,\det)$ respectively.

\begin{theorem}\label{thm-gt}
 The modular variety $M_{A, 2d}$ (or the modular variety $
 M^{(2d)}_{Q}$) is of general type if there exists a cusp form $F\in
 S_k(\Og(L_{A ,2d}), \det)$, (or $F\in S_k(\TOpl(L_{Q}), \det)$) of
 weight $k<21$ that vanishes of order at least one along the branch
 divisor of the modular projection
$$
\pi\colon \cD(L_{A, 2d}) \to \Og(L_{A ,2d})\setminus \cD(L_{A, 2d})
$$
(or the analogous projection for $\TOpl(L_{Q})$).
\end{theorem}

This is a particular case of Theorem 1.1 in \cite{GHS1}.

The dimension of the modular variety is smaller than $26$. Then we can
use the \emph{quasi pull-back} (see \cite{Bo}, \cite{BKPS}, \cite{Ko},
\cite{GHS1} and equation~\eqref{eq-pb} below) of the Borcherds modular
form
$$
\Phi_{12}\in M_{12}(\Orth^+(II_{2,26}),\det)\quad \text{where}\quad
II_{2,26}\cong 2U\oplus 3E_8(-1).
$$
We note that $\Phi_{12}(Z)=0$ if and only if there exists $r\in
II_{2,26}$ with $r^2=-2$ such that $(r,Z)=0$.  {Moreover, the
  multiplicity of the divisor of zeroes} of $\Phi_{12}$ is~$1$ (see
\cite{Bo}). We used the quasi pull-back of $\Phi_{12}$ in order to
construct cusp forms of small weight on the moduli spaces of polarised
$\Kthree$ surfaces (see \cite{GHS1}) and on moduli spaces of
split-polarised symplectic manifolds of $\Kthree^{[2]}$-type (see
\cite{GHS2}), which have dimension $19$ and $20$ respectively. The
present case is of dimension $21$. The non-split case is similar to
the cases considered in \cite{GHS1}--\cite{GHS2} (see also the example
at the end of this section) but the split case is different from the
previous ones because we need a cusp form with respect to the modular
group $\Og(L_{A ,2d})$, which is strictly larger than the stable
orthogonal group $\TOpl(L_{A, 2d})$. For this reason we will
concentrate in this paper on the split case.
\smallskip

Let $S\subset E_8(-1)$ be a sublattice (primitive or not) of rank $3$. For
our present purpose we take the sublattice of polarisations
$S=A_2(-1)\oplus \latt{ -2d}$ or $S=Q(-1)$ from Theorem
\ref{thm-O(L,h)}. The choice of $S$ in $E_8(-1)$ determines an
embedding of $L_{S}=2U\oplus 2E_8(-1)\oplus S$ into $II_{2,26}$. The
embedding of the lattice also gives us an embedding of the domain
$\cD(L_{S})\subset \PP (L_{S}\otimes \CC)$ into $\cD(II_{2,26})
\subset \PP(II_{2,26} \otimes \CC)$.

We put $R_S=\{r\in E_8(-1)\mid r^2=-2,\ (r, S)=0\}$, and $N_S=\#
R_S$. Then the quasi pull-back of $\Phi_{12}$ is given by the
following formula:
\begin{equation}\label{eq-pb}
 \left. F_S=
 \frac{\Phi_{12}(Z)}{
 \prod_{\{r\in R_S,\ r>0\}} (Z, r)}
 \ \right\vert_{\cD(L_{S})}
 \in M_{12+\frac{N_S}2}(\Tilde\Orth^+(L_{S}),\, \det).
\end{equation}
We fix a system of simple positive roots in $E_8(-1)$ and the notation
$r>0$ in the above formula means that we take the positive roots in
$R_S$, i.e.\ we pick only one root in any $A_1\subset R_S$. (The particular
choice of a system of the simple roots is not important.) The form
$F_S$ is a non-zero modular form of weight $12+\frac{N_S}2$. By
\cite[Theorems 6.2 and 4.2]{GHS1} it is a cusp form if $N_S\neq 0$,
since any isotropic subgroup of the discriminant form of the lattice
$L_{S}$ is cyclic, by Theorem~\ref{thm-O(L,h)}.
\smallskip

\noindent {\bf Example 2.} \emph{The smallest non-split
  polarisations.} We illustrate the method of Theorem~\ref{thm-gt}
together with the quasi pull-back construction for the polarisations
from Example~1 of \S \ref{sec3}. For the first three polarisations the
cusp form $F_Q$ is of weight $32$, $23$ and $22$ respectively. But for
the lattice $Q$ of determinant $22$ ($h^2=66$) we have a cusp form of
small weight $19<21$
$$
F^{(22)}_Q\in S_{19}(\Tilde\Orth^+(L_{Q}),\, \det).
$$
To apply Theorem \ref{thm-gt} we need a cusp form of small weight
with zero along the ramification divisor of the modular
projection. According to \cite[Corollary 2.13]{GHS1} this divisor is
determined by plus or minus reflections $\pm \sigma_r$ in the
corresponding modular group. If $\sigma_r$ is a reflection in this
group then $F^{(22)}_Q(\sigma_r(Z))=-F^{(22)}_Q(Z)$ and
$F^{(22)}(Z)=0$ if $(Z,r)=0$. If $-\sigma_r\in \Tilde\Orth^+(L_{Q})$
then $\det(-\sigma_r)=1$ because the dimension is odd. The weight of
$F^{(22)}_Q$ is also odd,
i.e.\ $F^{(22)}_Q(-Z)=-F^{(22)}_Q(Z)$. Therefore
$$
-F^{(22)}_Q(\sigma_r(Z))=F^{(22)}_Q(-\sigma_r(Z))=\det(-\sigma_r)F^{(22)}_Q(Z)
=F^{(22)}_Q(Z)
$$
and $F^{(22)}_Q$ vanishes along the divisor defined by $r$. Applying
Theorem \ref{thm-gt} we obtain

\begin{corollary}\label{prop-gt22}
The modular variety $M^{(66)}_{Q}$ is of general type.  Every
component $\cM_h^0$ of the moduli space $\cM_h$ of $10$-dimensional
polarised O'Grady varieties with non-split polarisation $h$ of
Beauville degree $h^2=66$ is of general type.
\end{corollary}

Any vector $l$ of length $12$, $30$ or $48$ with $\div(l)=3$ is
orthogonal to at least $20$ roots in $E_6$.  Hence we cannot apply the
low weight cusp form trick.  We conjecture that \emph{for the three
  lowest non-split polarisations, of Beauville degrees $2d=12$, $30$
  and $48$, the corresponding moduli spaces are unirational}. Using
the arithmetic and analytic methods developed in
\cite{GHS1}--\cite{GHS2} we hope to prove that for other non-split
polarisations the moduli spaces are of general type. In this paper we
study the split polarisation because this case is very different and
has new phenomena appearing.
\smallskip

The Weyl group of $E_8$ acts transitively on the sublattices
$A_2$. Let us fix a copy of $A_2(-1)$ in $E_8(-1)$. Then
$(A_2(-1))^\perp_{E_8(-1)}\cong E_6(-1)$. Let $l\in E_6(-1)$ satisfy
$l^2=-2d$. We denote the quasi pull-back $F_S$ for $S=A_2(-1)\oplus
\latt{l} $ by $F_l$. The problem is to find such a vector $l$ in
$E_6(-1)$ that yields a modular form with respect to the larger group
$\Og(L_{A ,2d})$.

\begin{lemma}\label{lem-Gmod}
Let us assume that $l\in E_6(-1)$, $l^2=-2d$, is invariant with
respect to the involution of the Dynkin diagram of $E_6(-1)$. Then the
quasi pull-back $F_l$ is modular with respect to $\Og(L_{A ,2d})$.
\end{lemma}
\begin{proof}
We see that $\Og(L_{A ,2d})=\latt{\TOpl(L_{A ,2d}), \sigma_6}$ where
$\sigma_6$ is a reflection with respect to any $-6$-vector in
$A_2(-1)$ (see \eqref{eq-Og}). The involution $\sigma_6\in W(G_2(-1))$
induces $-\id$ on the first component $D(A_2(-2))$ of the discriminant
group $D(L_{A,2d})$. The Weyl group $W(E_6)$ is a subgroup of
index~$2$ in $\Orth(E_6)$. The involution $J$ of the Dynkin diagram of
the fixed system of simple roots of $E_6(-1)$ induces $-\id$ on
$D(E_6(-1))$, which is also cyclic of order $3$.  Using the fact that
$(A_2)^\perp_{E_8}\cong E_6$ we can extend the element
$J_6=(\sigma_6,J)$ to an element in $\Orth(E_8) < \Orth^+(II_{2,26})$
where we consider $\sigma_6$ as an element in $\Orth^+(2U\oplus
2E_8(-1)\oplus A_2(-1))$. Let us introduce the coordinates $(Z_1,z_2,
Z_3)\in \cD(II_{2,26})$ corresponding to the sublattice
$$
\bigl(2U\oplus 2E_8(-1)\oplus A_2(-1)\bigr)\oplus \latt{l} \oplus
l^\perp_{E_6(-1)} \subset II_{2,26}
$$
where $z_2\in l\otimes \CC$ and $Z_3\in l^\perp_{E_6(-1)}\otimes
\CC$.  We calculate the function
$$
\left.\frac{\Phi_{12}(J_6(Z_1,z_2,Z_3))} {\prod_{\{r\in R_l,\,r>0\}}
(J_6(Z_1,z_2,Z_3), r)}
 \ \right\vert_{\cD(L_{A,2d})}
$$
where $R_l=\{r\in E_6(-1)\,|\,r^2=-2,\ (r,l)=0\}$ is the set of roots
in $E_8(-1)$ orthogonal to $S=A_2(-1)\oplus \latt{l}$.  First, we
find that it is equal to
$$
\left. \frac{\Phi_{12}((\sigma_6 Z_1,z_2, J(Z_3)))}
{\prod_{\{r\in R_l,\,r>0\}} ((\sigma_6 Z_1, z_2, J(Z_3)), r)}
 \ \right\vert_{\cD(L_{A,2d})}=F_l(\sigma_6 (Z_1,z_2))
$$
because $J(l)=l$ and $J_6(z_2)=z_2$.  Second, using the fact that
$\Phi_{12}$ has character $\det$ we find that the same function is
equal to
$$
\left.\frac{(\det J_6)\Phi_{12}(Z_1,z_2,Z_3)} {\prod_{\{r\in R_l,\,r>0\}}
((Z_1,z_2,Z_3), J(r))}
 \ \right\vert_{\cD(L_{A,2d})}=-F_l((Z_1,z_2))
$$
because $\det J=1$, $\det \sigma_6=-1$, $\det J_6=-1$ and the
involution $J$ permutes the positive roots in $l^\perp_{E_6}$.  We
note also that $(\sigma_6 Z_1,z_2, J(Z_3), r)=(J(Z_3),
r)_{E_6}=(Z_3,J(r))_{E_6}$.  Therefore
\begin{equation}\label{eq-Fl}
F_l\in S_{12+\frac{N_l}{2}}(\Og(L_{A ,2d}), \det)
\end{equation}
where $N_l=\#\{r\in E_6(-1)\,|\,r^2=-2,\ (r,l)=0\}$.
\end{proof}

The weight of $F_l$ is smaller than $21$ if $N_l<18$. In Section
\ref{sec4} we determine all $d$ for which there exists a
$(-2d)$-vector in $E_6(-1)$ invariant with respect to the automorphism
of the Dynkin diagram. In the next lemma we study the ramification
divisor of the modular projection of $\Og(L_{A ,2d})$. We studied this
divisor for the modular groups $\TOpl(L)$ in \cite[Proposition
  3.2]{GHS1} but the ramification divisor of $\Og(L_{A ,2d})$ is much
larger.
\begin{lemma}\label{lem-refl}
If $-\sigma_r\in \Og(L_{A ,2d})$, then $r^2=-2d\ $ and $\ \div(r)=2d$,
or $r^2=-6d$ and $\div(r)=3d$, or $r^2=-2d$ and $\div(r)=d$.
\end{lemma}
\begin{proof}
Let $r\in L_{A,2d}$ be a primitive vector and $r^2=-2e$. If
$\sigma_r\colon v\mapsto v-\frac{2(v,r)}{(r,r)}r\in \Orth^+ (\LA)$
then
$$
\div (r)\mid r^2 \mid 2\div (r)\quad{\rm and }\quad \div (r)\mid \lcm(3,2d).
$$

We assume that $-\sigma_r\in \Opl(L_{A,2d})$. Then
$\sigma_r|_{D(\latt{ -2d})}=-\id$ and for any $v\in L_{A,2d}^\vee$ we
have
$$
\sigma_r(v)+v=2v-\frac{2(v,r)}{(r,r)}r=2v-(v,r)\frac{r}e\in
A_2(-1)^\vee+L_{A,2d}
$$
where $(v,r)\in \ZZ$. This is true because we have no
$D(\latt{-2d})$-part in the sum $\sigma_r(v)+v$.  In particular, there
are the following relations between abelian groups
$$
2\cdot D(L_{A,2d})\cong \ZZ/3\ZZ\oplus \ZZ/d\ZZ\, <\, \ZZ/3\ZZ+\ZZ/e\ZZ,
$$
where the sum of the subgroup is taken in the discriminant
group. Therefore $d\vert e$. We have
$$
d\mid e\mid \div(r)\mid 2e\quad\text{ and }\quad \div(r)\mid \lcm(3,2d).
$$
Our aim is to calculate the two lattices
\begin{equation}
L_{A,2d}^{(r)}=r^\perp_{L_{A,2d}}\text{ and }
T_{r,d}=(L_{A,2d}^{(r)})^\perp_{II_{2,26}}.
\end{equation}
According to Lemma \ref{lem-Lh} we have
$$
\det T_{r,d}=\det L_{A,2d}^{(r)}= \frac{12de}{(\div(r))^2}.
$$
Analysing all possible $e$ and $\div(r)$ we see that $\det T_{r,d}$ is
a divisor of~$12$. The possible cases are
\begin{alignat*}{4}
e=\ d,&\quad r^2=\ 2d,&\quad \div(r)=\ d,&\quad \det T_{r,d}=12;\\
e=\ d,&\quad r^2=\ 2d,&\quad \div(r)=2d, &\quad \det T_{r,d}=\ 3;\\
e= 2d,&\quad r^2=\ 4d,&\quad \div(r)=2d, &\quad \det T_{r,d}=\ 6;\\
e= 3d,&\quad r^2=\ 6d,&\quad \div(r)=3d, &\quad \det T_{r,d}=\ 4;\\
e= 3d,&\quad r^2=\ 6d,&\quad \div(r)=6d, &\quad \det T_{r,d}=\ 1;\\
e= 6d,&\quad r^2=12d,&\quad \div(r)=6d, &\quad \det T_{r,d}=\ 2.
\end{alignat*}
In \cite[Table I]{CS} one can find all indecomposable lattices of
small rank and determinant. Analysing all lattices of determinant
$\det \mid 12$ and of rank $n\le 6$ we find the five classes
\begin{equation}
\det=3,\ E_6;\quad \det=4,\ D_6;\quad \det=12,\ A_5\oplus A_1,
\ D_4\oplus A_2,\ [D_5\oplus \latt{12}]_2
\end{equation}
where $[D_5\oplus \latt{12}]_2$ denotes an overlattice of order $2$ of
$D_5\oplus \latt{12}$. The root system of $[D_5\oplus \latt{12}]$ is
$D_5$. The formula for $\det T_{r,d}$ given above shows that only the
cases mentioned in the lemma are possible.
\end{proof}
\begin{corollary}\label{lem-branch0}
Let $l$ be as in Lemma \ref{lem-Gmod}. We assume that $N_l<18$. Then
the quasi pull-back $F_l$ vanishes along the ramification divisor of
the modular projection
$$
\pi\colon \cD(L_{A, 2d}) \to \Og(L_{A ,2d})\setminus \cD(L_{A, 2d}).
$$
\end{corollary}
\begin{proof}
The components of the branch divisor are
$$
\cD_r=\{\,[Z]\in \cD(L_{A,2d})\,|\, (r,Z)=0\,\}
$$
where $r\in L_{A,2d}$ and $\sigma_r$ or $-\sigma_r$ is in
$\Og(L_{A,2d})$ (see \cite[Corollary 2.13]{GHS1}). If $\sigma_r\in
\cD(L_{A, 2d})$, then $F_l$ vanishes along $\cD_r$ because $F_l$ is
modular with character $\det$. Let $-\sigma_r\in \Opl(L_{A,2d})$. The
divisor $\cD_r$ coincides with the homogeneous domain
$\cD(L_{A,2d}^{(r)})$. The Borcherds modular form $\Phi_{12}$ vanishes
of order $N/2$ where $N\ge |R(D_4\oplus A_2)|=30$ is the number of
roots in the lattice $\det T_{r,d}$. Since $N_l<18$ then the form
$F_l$ vanishes along $\cD_r$ with order at least $7$.
\end{proof}

\section{The $2d$-vectors in $E_6$ and the root system $F_4$}
\label{sec5}

In this section we finish the proof of Theorem~\ref{thm-main}.  To
prove it we use Theorem~\ref{thm-gt}, Lemma~\ref{lem-Gmod} and
Lemma~\ref{lem-branch0}. We want to know for which $2d>0$ there exists
a vector $l\in E_6$ of length $l^2=2d$, invariant with respect to the
involution $J$ of the Dynkin diagram of $E_6$ and orthogonal to at
least $2$ and at most $16$ roots in $E_6$. The answer is given in the
next theorem.

\begin{theorem}\label{thm-d}
A $J$-invariant vector $l$ of length $l^2=2d$ that is orthogonal to at
least $2$ and at most $16$ roots in $E_6$ exists if $d$ is not equal
to $2^n$ where $n\ge 0$.
\end{theorem}

We give the proof of the theorem in Lemmas~\ref{lem-E6}--\ref{lem-18}
below. We use the notation $A_n$, $D_n$ or $E_n$ both for a lattice
and for its root system because it is always clear from the context
which is meant. We consider the Coxeter basis of simple roots in the
lattice $E_6=\latt{\alpha_1,\dots, \alpha_6}$ (see \cite[Table V]{Bou})
\newpage

\noindent\hskip 2cm
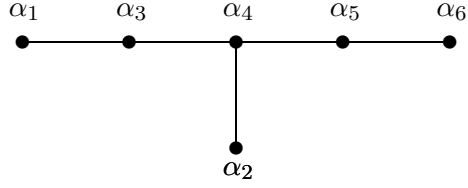
\begin{picture}(300,10)(55,10) \put(100,0){\circle*{5}}
\put(95,10){$\alpha_1$} \put(100,0){\vector(1,0){42}} \put(140,0){\circle*{5}}
\put(135,10){$\alpha_3$} \put(140,0){\vector(1,0){42}} \put(180,0){\circle*{5}}
\put(175,10){$\alpha_4$} \put(180,1){\vector(0,-1){43}}
\put(180,-40){\circle*{5}} \put(175,-50){$\alpha_2$} \put(175,-50){$\alpha_2$}
\put(180,0){\vector(1,0){42}} \put(220,0){\circle*{5}} \put(215,10){$\alpha_5$}
\put(220,0){\vector(1,0){42}} \put(260,0){\circle*{5}} \put(255,10){$\alpha_6$}
\end{picture}
\vskip2cm where
\begin{gather*}
\alpha_1=\frac 1{2}(e_1+e_8)-\frac 1{2}(e_2+e_3+e_4+e_5+e_6+e_7),\\
\alpha_2=e_1+e_2,\quad \alpha_k=e_{k-1}-e_{k-2}\ \ (3\le k\le 6)
\end{gather*}
and $(e_1,\dots, e_8)$ is a Euclidean basis in $\ZZ^8$. To get the
extended Dynkin diagram one has to add the maximal root
\begin{gather*}
\tilde \alpha=\frac 1{2}(e_1+e_2+e_3+e_4+e_5-e_6-e_7+e_8)\\
=\alpha_1+2\alpha_2+2\alpha_3+3\alpha_4+2\alpha_5+\alpha_6.
\end{gather*}
Then $(-\tilde\alpha, \alpha_2)=-1$ and $-\tilde\alpha$ is orthogonal
to all other simple roots.

In the Euclidean basis $(e_i)$ we have the following representation of
$E_6$
\begin{equation}\label{eq-E6}
E_6=\{l=x_1e_1+\dots+x_5e_5+x_6(e_6+e_7-e_8)\},
\end{equation}
$$
l^2=x_1^2+\dots+x_5^2+3x_6^3
$$
where the $x_i$ are either all integral or all half-integral, and in
both cases $x_1+\dots+x_6$ is an even integer. We recall that
$$
\Aut(E_6)=W(E_6)\times \Aut ({\rm Dynkin\ diagram\ of\ } E_6)
$$
where the second factor is the cyclic group of order $2$ generated by
the involution $J$ given by $J(\alpha_1)=\alpha_6$,
$J(\alpha_3)=\alpha_5$, $J(\alpha_4)=\alpha_4$,
$J(\alpha_2)=\alpha_2$.

\begin{lemma}\label{lem-E6}
The involution $J$ defines sublattices $E_6^{J,+}\oplus E_6^{J,-} \subset
E_6$ of index $4$ in $E_6$, where
\begin{gather*}
E_6^{J,+}=\{l\in E_6\,|\, J(l)=l\}\cong D_4,\\
E_6^{J,-}=\{l\in E_6\,|\, J(l)=-l\}\cong A_2(2)
\end{gather*}
and $A_2(2)$ is the lattice with the quadratic form
$\left(\smallmatrix 4&-2\\-2&\ 4\endsmallmatrix\right)$ (the
renormalisation of the lattice $A_2$ by $2$).
\end{lemma}
\begin{proof}
From the definition of $J$ we have $E_6^{J,+}=\latt{\alpha_2,
  \,\alpha_4,\, \alpha_1+\alpha_6,\, \alpha_3+\alpha_5}$. This has
another basis, namely
\begin{gather*}
E_6^{J,+}=\latt{\alpha_2,\, \alpha_4,\, \alpha_3+\alpha_4+\alpha_5,\,
(\alpha_1+\alpha_6)+2(\alpha_3+\alpha_4+\alpha_5)+2\alpha_2+\alpha_4
}\\
=\latt{\alpha_2, \,\alpha_4,\, \alpha_3+\alpha_4+\alpha_5,\, -\tilde\alpha
} \cong D_4
\end{gather*}
where $\alpha_2$ is the central root of the Dynkin diagram of
$D_4$. We denote $E_6^{J,+}$ by $D_4^+$.

If $J(u)=u$ and $J(v)=-v$ then $(u,v)=-(u,v)=0$. Therefore
$$
E_6^{J,-}=(D_4^+)^\perp_{E_6} \supseteq \latt{\alpha_1-\alpha_6,\,
  \alpha_3-\alpha_5} \cong A_2(2)=
\begin{pmatrix}\ 4&-2\\-2&\ 4\end{pmatrix}.
$$
A direct calculation shows that we have equality in the above
inclusion of lattices.  Then we have $\det D_4=4$ and $\det
A_2(2)=12$, so $[E_6,D_4^+\oplus A_2(2)]=4$.
\end{proof}

In what follows we need some properties of the root systems $D_4$ and
$F_4$.  The lattice $D_n$ is a sublattice of the Euclidean lattice
$\ZZ^n$
$$
D_n=\{l=(x_1,\dots,x_n)\in \ZZ^n\,|\, x_1+\dots+x_n\in 2\ZZ\}.
$$
The lattice $D_4$ contains the twenty-four $2$-roots
$$
R_2(D_4)=\{\pm (e_i\pm e_j),\ 1\le i<j\le 4\}
$$
which form the root system $D_4$. But the lattice $D_4$ contains also
the twenty-four $4$-roots
$$
R_4(D_4)=\{\pm e_1\pm e_2\pm e_3\pm e_4,\ \pm 2e_i,\ 1\le i\le 4\}.
$$
By definition of the root system $F_4$ equals
$$
F_4=R_2(D_4)\cup R_4(D_4).
$$
The Weyl group of $F_4$ coincides with the orthogonal group of the
lattice $D_4$:
$$
\Orth(D_4)=W(F),\qquad W(F_4)/W(D_4)\cong \Aut({\rm
  Dynkin\ diagram\ of\ } D_4)\cong S_3.
$$
\begin{lemma}\label{lem-J}
Let $J$ be the involution of the Dynkin diagram of $E_6$.
\newline
\noindent 1) For any root $r\in R_2(E_6)$ we have
$$
J(r)\ne r \Leftrightarrow (J(r),r)=0.
$$
\noindent
2) For $D_4^+=E_6^{J,+}$ we have
$$
R_4(D_4^+)=\{\,r+J(r)\,|\,r\in R_2(E_6),\ r\ne J(r)\}.
$$
\noindent 3) Let $l\in D_4^+$ be orthogonal to a vector $l_4\in
R_4(D_4^+)$.  Then $l$ is orthogonal to the roots $r$ and $J(r)$ from
$E_6$ such that $l_4=r+J(r)$ and $r\ne J(r)$.
\end{lemma}
\begin{proof}
1) Lemma \ref{lem-E6} gives us the following inclusion of lattices:
\begin{equation}\label{eq-EJ}
D_4^+\oplus A_2(2)\subset E_6\subset E_6^\vee \subset
(D_4^+)^\vee\oplus A_2(2)^\vee.
\end{equation}
We proved above that
$$
[E_6: (D_4^+\oplus A_2(2))]=[D_4^\vee : D_4]=\det D_4=4.
$$
It is easy to see that
$$
D_4^\vee/D_4=\{ 0,\ e_1+D_4,\ \frac 1{2}(e_1+e_2+e_3\pm
e_4)+D_4\}\cong \ZZ/2\ZZ \times \ZZ/2\ZZ
$$
where
$$
q_{D_4}(e_1+ D_4)= q_{D_4}\bigl(\frac 1{2}(e_1+e_2+e_3\pm e_4)+ D_4\bigr)
\equiv 1 \mod 2\ZZ.
$$
Analysing the discriminant form $A_2(2)^\vee/A_2(2)$ we see that it
contains only three classes $\frac{1}2 a$, $\frac{1}2 b$ and
$\frac{1}2 (a+b)$ modulo $A_2$ (where $a$, $b$ are simple roots in
$A_2$) of square $1 \mod 2\ZZ$.  Using (\ref{eq-proj-pH}) we see that
the natural projection $E_6/(D_4^+\oplus A_2(2))$ onto $D_4^\vee/D_4$
is surjective. It follows that if
$$
l\in E_6,\quad l=l_+^*+l_{-}^*,\quad \text{where}\quad l_+^*\in
(D_4^+)^\vee,\
l_{-}^*\in A_2(2)^\vee,\ l_{+}^*\not\in D_4^+
$$
then $(l_{+}^*,l_{+}^*)\equiv 1 \mod 2\ZZ$.

Let consider this representation $r_+^*+r_{-}^*$ for a root $r$ in
$E_6$. Then $r^2=(r_+^*)^2+(r_{-}^*)^2=2$ and the second component
$r_{-}^*$ is non-trivial if and only if $(r_+^*)^2=(r_{-}^*)^2=1$
according to the argument above. Then $J(r)\ne r$ if and only if $(r,
J(r))=(r_+^*)^2-(r_{-}^*)^2=0$.
\smallskip

2) We showed in Lemma \ref{lem-E6} that $E_6$ contains exactly $24$
$J$-invariant roots of $D_4^+$. Therefore there are $72-24=48$
non-invariant roots. For any non-invariant root $r$ we proved in 1)
that $(r,J(r))=0$. This gives us $24$ pairs $(r,J(r))$ of
non-invariant roots satisfying $(r+J(r))^2=4$ and $r+J(r)\in
D_4^+$. To show that there is a bijection between the $J$-pairs and
$4$-roots in $D_4^+$ one can simply pick $\alpha_1+J(\alpha_1)$ and
take into account the fact that the Weyl group of $D_4$ acts
transitively on the set of $4$-vectors in $D_4$.
\smallskip

3) If $l\in D_4^+$ then $(l,r)=(l,J(r))$ for any root. Therefore
$2(l,r)=(l, r+J(r))=0$.
\end{proof}

\begin{lemma}\label{lem-1root}
For any positive integer $d$ there exists a vector $l_{2d}\in
D_4^+=E_6^{J,+}$ of square $2d$ which is orthogonal to at least one
root in $E_6$.
\end{lemma}
\begin{proof}
We denote by $N_L(2d)$ the number of vectors of square $2d$ in a
positive definite lattice $L$. We consider two cases: a vector
$l_{2d}$ is orthogonal to a $J$-invariant root $r_{J}$ or to a
non-$J$-invariant root $r_{n}$. In the first case $l_{2d}\in
(r_J)^\perp_{D_4^+}\cong 3A_1$. (See the fourth case in the proof of
Lemma \ref{lem-18} below.) Then
$$
N_{3A_1}(2d)=r_3(d)
$$
where $r_3(d)$ is equal to the number of representations of $d$ as a
sum of three squares. It is classically known that
\begin{equation}\label{eq-3sq}
r_3(4^md)=r_3(d) \quad{\rm and }\quad r_3(d)>0 \quad {\rm if } \ d\ne
2^{2m}(8n+7).
\end{equation}
If $(l_{2d}, r_n)=0$ then $(l_{2d}, r_n+J(r_n))=0$ where
$r_n+J(r_n)=l_4\in D_4^+$. But
$$
(l_4)^\perp_{D_4^+}\cong A_3.
$$
This follows from the form of the extended Dynkin diagram of
$D_4$. For $l_4$ we can take the alternating sum of two orthogonal
simple roots. Then the three other roots of the extended diagram form
the orthogonal complement of $l_4$. We have $A_3\cong D_3$. According
to the definition of $D_3$ we have that $N_{A_3}(2d)=r_3(2d)$. The
last number is not zero if $d\ne 2^{2m-1}(8n+7)$. This and formula
(\ref{eq-3sq}) shows that for any $d$ we have
$N_{3A_1}(2d)+N_{A_3}(2d)>0$. This proves the lemma.
\end{proof}

\begin{lemma}\label{lem-18}
Let $l_{2d}$ be a vector as in Lemma~\ref{lem-1root}. Then the number
of roots in $E_6$ orthogonal to $l_{2d}$ is smaller than $18$ if and
only if $d$ is not equal to $2^n$ where $n\ge 0$.
\end{lemma}
\begin{proof}
Let us assume that $|R_2((l_{2d})^\perp_{E_6})|\ge 18$. The root
systems of rank at most $5$ having at least $18$ roots are
$$
A_5,\ D_5,\ A_4\oplus A_1,\ D_4\oplus A_1,\ A_3\oplus A_2,
\ A_4,\ D_4.
$$

1) The cases of $A_3\oplus A_2$ and $D_4\oplus A_1$ are not
possible. $W(E_6)$ acts transitively on the roots and on the
$A_2$-sublattices of $E_6$. We have $(A_1)^\perp_{E_6}\cong A_5$ and
$(A_2)^\perp_{E_6}\cong A_2\oplus A_2$. But $A_5$ does not contain
$D_4$ and $A_2\oplus A_2$ does not contain $A_3$.
\smallskip

2) Let us assume that $R_2((l_{2d})^\perp_{E_6})=A_4$ or $A_4\oplus
A_1$.  We show that neither case is possible. The vector $l_{2d}$ is
$J$-invariant.  Therefore $J(A_4)=A_4$. The lattice $A_4$ is generated
by its simple roots $a_1$, $a_2$, $a_3$ and
$a_4$:
\vskip.5cm \noindent\hskip 1cm
\begin{picture}(300,10)(55,10)
\put(120,10){\circle*{5}} \put(115,0){$a_1$} \put(120,10){\vector(1,0){62}}
\put(180,10){\circle*{5}} \put(170,0){$a_2$} \put(180,10){\vector(1,0){62}}
\put(240,10){\circle*{5}} \put(225,0){$a_3$} \put(240,10){\vector(1,0){62}}
\put(300,10){\circle*{5}} \put(295,0){$a_4$}
\end{picture}
\vskip 1cm
First we assume $a_1\ne J(a_1)$ and $J(a_4)\ne a_4$. Then
$(a_1,J(a_1))=(a_4,J(a_4))=0$ according to Lemma
\ref{lem-J}. Therefore we have $J(a_4)\in \latt{a_1, a_2}$ and
$J(a_1)\in \latt{a_3, a_4}$. If $J(a_1)\ne \pm a_4$ then $A_4$
contains two orthogonal sublattices $\latt{ a_1, J(a_4)}$ and
$\latt{a_4, J(a_1)}$ isomorphic to $A_2$, which is impossible.

If $J(a_1)=\pm a_4$ then $0=(J(a_1),J(a_3))=(\pm a_4,J(a_3))$ and
$J(a_3)\in \latt{a_1, a_2}$. But $J(a_3)\ne \pm a_1$ and we obtain
that $J(a_3)\ne a_3$ and $(J(a_3), a_3)\ne 0$. This contradicts Lemma
\ref{lem-J}.  Therefore we can assume that $a_1= J(a_1)$ or $a_4=
J(a_4)$. If $a_1= J(a_1)$ then $(a_1, J(a_4))=0$ and $J(a_4)\in
\latt{a_3, a_4}$. It follows that $J(a_4)=a_4$. An analogous argument
shows that $J(a_3)=a_3$ and $J(a_2)=a_2$.  Therefore $J$ is the
identity on $A_4$ and we obtain that $A_4$ is a sublattice of
$D_4^+=E_6^{J,+}$, which is impossible.  If
$R_2((l_{2d})^\perp_{E_6})=A_4\oplus A_1$ then again we have that
$J(A_4)=A_4$ and $J|_{A_4}=\id$.
\smallskip

3) We have mentioned above that $(A_1)^\perp_{E_6}\cong A_5$ and that
there is only one $W(E_6)$-orbit of $A_1$ in $E_6$.  Therefore
$(A_5)^\perp_{E_6}\cong A_1=\latt{2}$. Any non-zero vector $l\in A_1$
($l^2=2m^2$) will have the same orthogonal complement. Let us take a
$J$-invariant vector $l\in 3A_1$ such that $l^2=2^{2n+1}k^2$ where $k$
is odd. Then $N_{3A_1}(2)=r_3(1)=6$ and
$$
N_{3A_1}(2^{2n+1}k^2)=r_3(k^2)=
\sum_{f|k}r_3^{pr}({k^2}/{f^2})=r_3^{pr}(1)+\dots +r_3^{pr}(k^2),
$$
which is $>6$ if and only if $k>1$. Here we denote by $r_3^{pr}(n)$
the number of primitive representation of $n$ by three
squares. According to Gauss $r_3^{pr}(n)=0$ if and only if $n\equiv
0\mod 4$ or $n\equiv 7\mod 8$. Therefore if $2d=2^{2n+1}$ then any
$2d$-vector in $3A_1$ is a multiple of a root. If $l_{2d}\in A_3$ the
situation is quite similar. We conclude that for $2d=2^{2n+1}k^2$
there is a $2d$-vector which satisfies the conditions of the lemma if
and only if $k>1$.
\smallskip

4) We can compare the case when $l^{\perp}_{E_6}=D_5$ with the case of
$A_5$.  We have $(D_5)^\perp_{E_6} \cong\latt{12}$. To see this we
consider $(D_5)^\perp_{E_8}=A_3$ and $(A_2)^\perp_{E_8}=E_6$. There is
only one $W(A_3)$-orbit of $A_2$ in $A_3$ and $(A_2)^\perp_{A_3}\cong
\latt{ 12}$. This gives us the sublattice $A_2\oplus D_5\oplus
\latt{12}$ in $E_8$. But we can find another orbit of $12$-vectors in
$E_6$ by taking a copy of $A_2$ in $D_5$. In fact, the $12$-vector
corresponding to the decomposition $\latt{12}\oplus D_5 \subset E_6$ is not
$J$-invariant.  To get a $J$-invariant vector we take
$$
l_{12}^+=2\alpha_2+\alpha_4=2e_1+e_2+e_3\in E_6^{J,+}
$$
(see the diagram of $E_6$ above). The roots of $E_6$ are the vectors
$$
\pm e_i\pm e_j\ (1\le i< j\le 5),\quad \pm \frac{1}2\bigl(
e_8-e_7-e_6\pm e_1\pm\dots\pm e_5\bigr)
$$
where the number of minus signs in the last case is even. We see that
there are six integral and eight half-integral roots orthogonal to $l_{12}^+$.
Up to sign they are
$$
e_3-e_2,\ e_4- e_5,\ e_4+e_5;
$$
$$
\frac {1}{2}\bigl( e_8-e_7-e_6+e_1-e_2-e_3\pm (e_4+e_5)\bigr),\
\frac {1}{2}\bigl( e_8-e_7-e_6-e_1+e_2+e_3\pm (e_4-e_5)\bigr).
$$
These roots form a root system $A_1\oplus A_3$ where
$A_1=\latt{\alpha_4}= \latt{e_3-e_2}$ and
$$
A_3
=\latt{e_4-e_5,\ e_4+e_5,\ \frac {1}{2}\bigl(
e_8-e_7-e_6+e_1-e_2-e_3-e_4-e_5\bigr)}.
$$
Therefore in the case $2d=12$ a vector giving a low weight cusp form
does exist.
\smallskip

5) Let us assume that $R_2((l_{2d})^\perp_{E_6})=D_4$. Then
$J(D_4)=D_4$. We can fix a system of simple roots $(a_1,a_2,a_3,a_4)$
of $D_4$ ($a_2$ is the central root of the diagram).

First we prove that $J(a_2)=a_2$. Consideration of the extended Dynkin
diagram of $D_4$ shows that $(A_1)^\perp_{D_4}\cong 3A_1$. The four
pairwise orthogonal copies of $A_1$ in $D_4$ correspond to the
vertices of the extended Dynkin diagram of $D_4$: $a_1$, $a_3$, $a_4$
and $-\tilde a$ where $\tilde a=a_1+2a_2+a_3+a_4$ is the maximal root
of $D_4$ (see \cite[Table IV]{Bou}). If $J(b)\ne b$ for a root $b$
then $J(b)$ is orthogonal to $b$ (Lemma \ref{lem-J}). Therefore $J$
permutes the roots $a_1$, $a_3$, $a_4$ and $-\tilde a$ with some
possible changes of signs. Therefore
$$
J(2a_2)=J(\tilde a-a_1-a_3-a_4)=\pm (a_1+2a_2+a_3+a_4 \pm a_1\pm a_3\pm a_4)
$$
where all $\pm$ are independent.  The maximal root $\tilde a$ is the
only root re\-presented by a linear combination of the simple roots
having a coefficient greater than $1$. That leaves only two
possibilities: $J(2a_2)=\pm 2(a_1+a_2+a_3+a_4)$ or $J(2a_2)=\pm
2a_2$. The first of those two does not occur because the root
$a_1+a_2+a_3+a_4$ is not orthogonal to $a_2$. Therefore $J(a_2)=a_2$.

Let us assume that $J$ does not fix any of the four pairwise
orthogonal copies $A_1$ in $D_4$. Let $J(a_1)\ne \pm a_3$ (the other
cases are similar). Then the root $J(a_1+a_2+a_3)=a_2+J(a_1)+J(a_3)$
is not equal to the root $a_1+a_2+a_3$ and it is not orthogonal to it.
This contradicts Lemma \ref{lem-J}-3).  Therefore $J$ fixes at least
one $A_1$ among the four copies of $A_1$. So $J$ fixes at least two
copies, which form together with $a_2$ a root system $A_3$ on which
$J$ acts trivially. Therefore we have proved that if $l_{2d}\in E_6$,
$J(l_{2d})=l_{2d}$ and $R_2((l_{2d})^\perp_{E_6})=D_4$, then the
orthogonal complement of $l_{2d}$ in $D_4^+=E_6^{J,+}$ contains
$A_3$. But $(A_3)^\perp_{D_4}\cong \latt{4}$. To see this one bears in
mind two facts: $W(F_4)=\Orth (D_4)$ acts transitively on the set of
$4$-vectors in $D_4$ and
$$
\latt{a_3-a_4}^\perp_{D_4}=\latt{a_1,\,a_2,\,-\tilde a}\cong A_3.
$$
It follows that the vector $l_{2d}$ is a multiple of a $4$-vector
$l_4$ in $D_4^+$
$$
l_{2d}=ml_4,\quad l_4\in 3A_1\subset D_4^+\ {\rm or }\ l_4\in A_3\subset  D_4^+
$$
(see Lemma \ref{lem-1root}).

If $2d=4m^2$ then any $2d$-vector in $3A_1\subset D_4^+$ or in $A_3\subset D_4^+$ is
a multiple of a corresponding $4$-vector if and only if $2d=4\cdot
2^{2n}$. We use an argument similar to the case $d=1$ (see part 3) of the proof
above). If $2d=4\cdot 2^{2n}k^2$, with $k$ odd, then
$$
N_{3A_1}(4\cdot 2^{2n}k^2)=r_3(2k^2)=
\sum_{f|k}r_3^{pr}(2\frac{k^2}{f^2})=r_3^{pr}(2)+\dots +r_3^{pr}(k^2),
$$
which is $>r_3(2)=12$ if and only if $k>1$. This finishes the proof
of~Lemma \ref{lem-18} and of Theorem~\ref{thm-d}.
\end{proof}

We note that by a remark of Freitag \cite[Hilfssatz 2.1, Kap.~3]{Fr} one can calculate the
geometric genus of a modular variety using cusp forms of canonical
weight. In particular we have
$$
p_g(M_{A,2d})= \dim S_{21}(\Og(L_{A,2d}),\det).
$$
In the cases of polarised $\Kthree$ surfaces or polarised symplectic
varieties of type $\Kthree^{[2]}$ we constructed canonical
differential forms on the corresponding modular varieties using the
quasi-pullback of $\Phi_{12}$.  In the case considered in this paper
this is not possible. From the proof of Lemma~\ref{lem-18} we obtain
\begin{corollary}
1. There are no $J$-invariant $2d$-vectors in $E_6$ which are
orthogonal to exactly $18$ roots in $E_6$.
\newline
\noindent 2. There are no $\Og(L_{A ,2d})$-modular quasi-pullbacks of
$\Phi_{12}$ of weight $21$.
\end{corollary}
We think that cusp forms of canonical weight exist for $\Og(L_{A ,2d})$,
but we expect the Beauville degree of the polarisation to be rather large.
To prove that the modular variety $M_{A,2d}$ with  $d=2^n$ is of general
type for $n$ large we could  use the explicit formula for
the Mumford-Hirzebruch volume found in \cite{GHS3}. We conjecture that this
variety is not of general type for small $n$, for example,  for $n=0$, $1$,
$2$. An argument for this is given in Proposition \ref{pr-tilde} below.

The modular variety of symplectic $10$-dimensional O'Grady varieties
with a split polarisation is a $2:1$ quotient of the modular variety
$$
\TOpl(L_{A,2d})\setminus \cD(L_{A, 2d}) \to \Og(L_{A ,2d})\setminus \cD(L_{A,
2d})=M_{A,2d}
$$
because $[\Og(L_{A ,2d}): \TOpl(L_{A,2d})]=2$.
\begin{proposition}\label{pr-tilde}
The modular variety $\TOpl(L_{A,2d})\setminus \cD(L_{A, 2d})$ is of
general type if $d\not\in \{\,1,\,2,\,4\,\}$.
\end{proposition}
\begin{proof}
We only have to consider the series $2d=2^n$.  If $2d=2$, $4$ or $8$
then any vector $l$ of length $l^2=2d$ is orthogonal to at least $20$
roots. We have seen this for $2d=2$ and $2d=4$. The argument for
$2d=8$ is similar. Hence we cannot apply the low weight cusp form
trick here.

The lattice $L_{A,2d}$ for $2d=2^n$ with $n>5$ can be considered as a
sublattice of $L_{A,16}$, if $n$ is even, or of $L_{A,32}$, if $n$ is
odd.  Therefore the corresponding modular variety is a covering of
finite order of one of the two varieties for $2d=16$ or $32$. Hence it
is enough to prove that $\TOpl(L_{A,16})\setminus \cD(L_{A, 16})$ and
$\TOpl(L_{A,32})\setminus \cD(L_{A, 32})$ are of general type.

\noindent
1) Let $2d=16$. Using the representation \eqref{eq-E6} of $E_6$ we put
$l_{16}=3e_1+2e_2+e_3+e_4+e_5\in E_6$.  Inspection shows that there
are $12$ orthogonal roots ($6$ copies of $A_1$).  Three ``integral"
copies are
$$
e_3-e_4,\ e_4-e_5,\ e_3-e_5.
$$
Three ``half-integral" copies are $\frac{1}{2}\bigl(-e_1+e_2\pm
(e_3-e_4)+e_5-e_6-e_7+e_8\bigr)$ and
$\frac{1}{2}\bigl(-e_1+e_2+e_3+e_4-e_5-e_6-e_7+e_8\bigr)$.  Then
$(l_{16})^\perp_{E_6}\cong A_3$ where
$$
A_3= \langle\, \frac{1}{2}\bigl(-e_1+e+2-e_3+e_4+e_5-e_6-e_7+e_8\bigr),\
e_3-e_4,\ e_4-e_5\rangle.
$$
2) Let $2d=32$. We put $l_{32}=4e_1+3e_2+2e_3+e_6+e_7-e_8\in
E_6$. Then $(l_{32})^\perp_{E_6}\cong A_2\oplus A_1$ where
$A_1=\latt{e_4+e_5}$ and
$$
A_2= \langle\, \frac{1}{2}\bigl(e_1-e+2+e_3-e_4+e_5-e_6-e_7+e_8\bigr),\ e_4-e_5
\rangle.
$$
The quasi pull-backs of $\Phi_{12}$ to $2U\oplus 2E_8(-1)\oplus
A_2(-1)\oplus \latt{-2d}$ for the vectors $l_{16}$ and $l_{32}$) are
cusp forms of weights $18$ and $16$ respectively, for the groups
$\TOpl(L_{A,16})$ and $\TOpl(L_{A,32}))$). The set of plus or minus
reflections in $\TOpl(L_{A,2d})$ is a subset of the reflections
considered in Lemma~\ref{lem-refl}. Therefore we can prove that
$F_{l_{16}}$ (resp. $F_{l_{36}}$) vanishes on the branch divisor of
the modular projection using the arguments of the proof of
Corollary~\ref{lem-branch0}. To finish the proof we apply
Theorem~\ref{thm-gt}.
\end{proof}
\bibliographystyle{alpha}

\bigskip
\noindent
V.A.~Gritsenko\\
Universit\'e Lille 1\\
Laboratoire Paul Painlev\'e\\
F-59655 Villeneuve d'Ascq, Cedex\\
France\\
{\tt valery.gritsenko@math.univ-lille1.fr}
\bigskip

\noindent
K.~Hulek\\
Institut f\"ur Algebraische Geometrie\\
Leibniz Universit\"at Hannover\\
D-30060 Hannover\\
Germany\\
{\tt hulek@math.uni-hannover.de}
\bigskip

\noindent
G.K.~Sankaran\\
Department of Mathematical Sciences\\
University of Bath\\
Bath BA2 7AY\\
England\\
{\tt gks@maths.bath.ac.uk}
\end{document}